\documentclass[12pt,a4paper,reqno]{amsart}
\usepackage{ae,aecompl}
\pdfoutput=1
\usepackage[english]{babel}

\title[Fredholm modules for quantum projective spaces]{Bounded and unbounded Fredholm modules\\[10pt] for quantum projective spaces}
\date{20 March 2009; v2: 25 January 2010}

\author[Francesco D'Andrea]{Francesco D'Andrea} 
\address{\flushleft D{\'e}partement de Math{\'e}matique, Universit\'e Catholique de Louvain,
Chemin du Cyclotron 2, B-1348, Louvain-La-Neuve, Belgium}
\email{francesco.dandrea@uclouvain.be}
\author[Giovanni Landi]{Giovanni Landi \\[20pt]}
\address{\flushleft Dipartimento di Matematica e
Informatica, Universit\`{a} di Trieste, Via A. Valerio 12/1, I-34127
Trieste, Italy, and INFN, Sezione di Trieste, Trieste, Italy}
\email{landi@univ.trieste.it}

\newtheorem{prop}{Proposition}
\newtheorem{lemma}[prop]{Lemma}
\newtheorem{df}[prop]{Definition}

\newcommand{\arxiv}[1]{[arxiv:#1]}

\newcommand{\U}{\mathrm{U}}
\newcommand{\A}{\mathcal{A}}
\newcommand{\B}{\mathcal{B}}
\newcommand{\HH}{\mathcal{H}}
\newcommand{\N}{\mathbb{N}}
\newcommand{\Z}{\mathbb{Z}}
\newcommand{\R}{\mathbb{R}}
\newcommand{\C}{\mathbb{C}}
\newcommand{\CP}{\mathbb{C}\mathrm{P}}
\newcommand{\inner}[1]{\left<#1\right>}
\newcommand{\ket}[1]{\left|#1\right>}
\newcommand{\tr}{\mathrm{Tr}}
\newcommand{\ma}[2]{\textrm{\footnotesize
$\left(\!\rule{0pt}{#1}\smash[b]{\smash[t]{\begin{array}{cc}#2\end{array}}}\!\right)$}}

\begin{document}

\begin{abstract}
We construct explicit generators of the K-theory and K-homology of the coordinate algebra of `functions' on quantum projective spaces. We also sketch a construction of unbounded Fredholm modules, that is to say Dirac-like operators and spectral triples of any positive real dimension.
\end{abstract}

\maketitle

\tableofcontents

% ----------------------------------------------------------------------
\section{Introduction}\label{se:intro}

For quantum projective spaces $\CP^n_q$ we generalize some ideas of \cite{DL09} and give `polynomial' generators of its  K-theory,  that is 
projections whose entries are in the corresponding coordinate algebra $\A(\CP^n_q)$. Dually, we give generators of its K-homology via even Fredholm modules 
$(\A(\CP^n_q),\HH_{(k)},\gamma_{(k)},F_{(k)})$, for  $k=0,\ldots,n$. The `top' Fredholm module -- the only one for which the representation of $\A(\CP^n_q)$ is faithful -- can be realized as the `conformal class' of a spectral triple $(\A(\CP^n_q),\HH_{(n)}, \gamma_{(n)},D)$, namely we realize $F_{(n)}:=D|D|^{-1}$ as the `sign' of a Dirac-like operator. This procedure allows us to construct spectral triples of any summability $d\in\R^+$.

In the following, without loss of generalities, the real deformation parameter is restricted to be $0<q<1$. Also, by $*$-algebra we shall mean a unital involutive associative $\C$-algebra, and by representation of a $*$-algebra we always mean a unital $*$-representation.

The `ambient' algebra  for the quantum projective space is the coordinate algebra  $\A(S^{2n+1}_q)$ of the unit quantum spheres. This $*$-algebra is generated by $2n+2$ elements $\{z_i,z_i^*\}_{i=0,\ldots,n}$ with relations \cite{VS91}:
\begin{subequations}
\begin{align}
z_iz_j &=q^{-1}z_jz_i &&\forall\;0\leq i<j\leq n \;,\\
z_i^*z_j &=qz_jz_i^* &&\forall\;i\neq j \;,\\
[z_i^*,z_i] &=(1-q^2)\sum\nolimits_{j=i+1}^n z_jz_j^* 
    &&\forall\;i=0,\ldots,n-1 \;,\label{eq:3rd}\\
[z_n^*,z_n] &=0 \;, \qquad \; \\
z_0z_0^*+z_1z_1^* &+\ldots+z_nz_n^*=1 \;.
\end{align}
\end{subequations}

The $*$-subalgebra generated by $p_{ij}:=z_i^*z_j$ will be denoted $\A(\CP^n_q)$, and identified with the algebra of `polynomial functions' on the quantum projective space $\CP^n_q$. 
The algebra $\A(\CP^n_q)$ is made of (co)invariant elements for the $\U(1)$ (co)action 
$ z_i  \to \lambda z_i $ for $\lambda \in\U(1)$.
{}From the relations of $\A(S^{2n+1}_q)$ one gets analogous relations for 
$\A(\CP^n_q)$:
\begin{align*}
p_{ij}p_{kl}&=q^{\mathrm{sign}(k-i)+\mathrm{sign}(j-l)}\, p_{kl}p_{ij}
 &\hspace{-1.5cm}\mathrm{if}\;i\neq l\;\mathrm{and}\;j\neq k\;,\\
p_{ij}p_{jk}&=q^{\mathrm{sign}(j-i)+\mathrm{sign}(j-k)+1}\, p_{jk}p_{ij}-(1-q^2)
\textstyle{\sum_{l>j}}\, p_{il}p_{lk} &\mathrm{if}\;i\neq k\;,\\
p_{ij}p_{ji}&=
q^{2\mathrm{sign}(j-i)}p_{ji}p_{ij}+
(1-q^2)\left(\textstyle{\sum_{l>i}}\,q^{2\mathrm{sign}(j-i)}p_{jl}p_{lj}
-\textstyle{\sum_{l>j}}\,p_{il}p_{li}\right)
 &\mathrm{if}\;i\neq j\;,
\end{align*}
with $\mathrm{sign}(0):=0$.
The elements $p_{ij}$ are  the matrix entries of a projection $P= (p_{ij})$, that is $P^2=P=P^*$ or 
$\sum_{j=0}^n p_{ij} p_{jk} = p_{ik}$ and $p_{ij}^*=p_{ji}$. 
This projection has  $q$-trace:
\begin{equation}\label{q-tr}
\tr_q(P):=\sum\nolimits_{i=0}^n \, q^{2i} p_{ii}=1.
\end{equation}

The original notations of \cite{VS91} are obtained by setting $q=e^{h/2}$;
the generators of \cite{DD09} correspond to the replacement $z_i\to z_{n+1-i}$,
while the generators $x_i$ used in~\cite{HL04} are related to ours by
$x_i=z_{n+1-i}^*$ and by the replacement $q\to q^{-1}$.

Generators for the K-theory and K-homology of the spheres $S^{2n+1}_q$ are in \cite{HL04}; unfortunately, there is no canonical way to obtain generators of the K-theory and/or K-homology of subalgebras, unless they are dense subalgebras (for a pair $C^*$-algebra/pre-$C^*$-algebra, the $K$-groups coincide). That $S^{2n+1}_q$ and $\CP^n_q$ are truly different can be seen from the fact that the $K$-groups of the odd-dimensional spheres are $(\Z,\Z)$, regardless of the dimension, while for $\CP^n_q$ they are $(\Z^{n+1},0)$.

That $K_0(C(\CP^n_q))\simeq\Z^{n+1}$ can be proved by viewing the corresponding $C^*$-algebra $C(\CP^n_q)$, the universal $C^*$-algebra of $\CP^n_q$, as the Cuntz--Krieger algebra of a graph \cite{HS02}. The group $K_0$ is given as the cokernel of the incidence matrix canonically associated with the graph. The dual result for K-homology is obtained using the same techniques: the group $K^0$ is now the kernel of the transposed matrix \cite{Cun84}; this leads to 
$K^0(C(\CP^n_q))=\Z^{n+1}$.

In \cite{HS02}, somewhat implicitly, there appear generators of the $K_0$ groups of $C(\CP^n_q)$ as projections in $C(\CP^n_q)$ itself. Here, we give generators of $K_0(C(\CP^n_q))$ in the form of `polynomial functions', so they represent elements of $K_0(\A(\CP^n_q))$ as well. They are also equivariant, i.e.~representative of elements in $K_0^{\mathcal{U}_q(\mathfrak{su}(n+1))}(\A(\CP^n_q))$. Besides, we give $n+1$ Fredholm modules that are generators of the homology group $K^0(C(\CP^n_q))$.

% ----------------------------------------------------------------------
\section{K-homology}

A useful (unital $*$-algebra) morphism $\A(S^{2n+1}_q)\to\A(S^{2n-1}_q)$ given by the map $z_n\mapsto 0$, restricts to a morphism $\A(\CP^n_q)\to\A(\CP^{n-1}_q)$ and is heavily used in what follows. We also stress that representations of $\A(S^{2n+1}_q)$, with $z_i$ in the kernel for all $i>k$, and for a fixed $0\leq k<n$, are the pullback of representations of $\A(S^{2k+1}_q)$.

Here we seek Fredholm modules (and then irreducible $*$-representations) for $\CP^n_q$ that are not the pullback of Fredholm modules (resp. irreps) for $\CP^{n-1}_q$. So, we look for irreps of $S^{2n+1}_q$ for which $z_n$ is not in the kernel. These are given in \cite{HL04}, and classified by a phase: in particular, they are inequivalent as representation of $S^{2n+1}_q$ but give the same representation of $\CP^n_q$ (whatever the value of the phase is). In fact, we start with representations that are not exactly irreducible, but close to it: they are the direct sum of an irreducible representation with copies of the trival $a\mapsto 0$ representation.

\begin{df}
We use the following multi-index notation. We let $\underline{m}=(m_1,\ldots,m_n)\in\N^n$ and,
for $0\leq i<k\leq n$, we denote by $\underline{\varepsilon}^k_i\in\{0,1\}^n$ the array
$$
\underline{\varepsilon}_i^k:=
(\,\,\stackrel{i\;\mathrm{times}}{\overbrace{0,0,\ldots,0}}\,,
\stackrel{k-i\;\mathrm{times}}{\overbrace{1,1,\ldots,1}}\,,
\stackrel{n-k\;\mathrm{times}}{\overbrace{0,0,\ldots,0}}) \;.
$$
Let $\HH_n:=\ell^2(\N^n)$, with orthonormal basis $\ket{\underline{m}}$.
For any $0\leq k\leq n$, a representation $\pi^{(n)}_k:\A(S^{2n+1}_q)\to\B(\HH_n)$,  
is defined as follows (all the representations are on the same Hilbert space).
We set $\pi^{(n)}_k(z_i)=0$ for all $i>k\geq 1$, while for the remaining generators
\begin{align*}
%\pi^{(n)}_k(z_0)\ket{\underline{m}} &=\sqrt{1-q^{2(m_1+1)}}\ket{\underline{m}+\underline{\varepsilon}_0^k} \;,\\
\pi^{(n)}_k(z_i)\ket{\underline{m}} &=q^{m_i}\sqrt{1-q^{2(m_{i+1}-m_i+1)}}
   \ket{\underline{m}+\underline{\varepsilon}_i^k}\;, && \forall\;0\leq i\leq k-1 \;,\\
\pi^{(n)}_k(z_k)\ket{\underline{m}} &=q^{m_k}\ket{\underline{m}} \;,
\end{align*}
on the subspace $\mathcal{V}^n_k\subset\HH_n$ --- linear span of basis vectors 
$\ket{\underline{m}}$ satisfying the restrictions
\begin{equation}\label{eq:mconstr}
0\leq m_1\leq m_2\leq\ldots\leq m_k\;,\qquad m_{k+1}>m_{k+2}>\ldots>m_n\geq 0 \;,
\end{equation}
with $m_0:=0$ ---, and they are zero on the orthogonal subspace. When $k=0$, we define $\pi^{(n)}_0(z_i)=0$ if $i>0$, $\pi^{(n)}_0(z_0)\ket{\underline{m}}=\ket{\underline{m}}$ if $m_1>m_2>\ldots>m_n\geq 0$,
and $\pi^{(n)}_0(z_0)\ket{\underline{m}} =0$ in all the other cases.
\end{df}

As a $*$-algebra, $\A(\CP^n_q)$ is generated by the elements $p_{ij}$, with $i\leq j$ since $p_{ji}=p_{ij}^*$; in fact, from the tracial relation \eqref{q-tr} one of the generators on the diagonal, say $p_{nn}$, is redundant. The computation of $\pi^{(n)}_k(p_{ij})$ for $i\leq j$ gives the following.
If $\underline{m}$ satisfies (\ref{eq:mconstr}):
\begin{align*}
\pi^{(n)}_k(p_{ij})\ket{\underline{m}} &=
q^{m_i+m_j}\sqrt{1-q^{2(m_{i+1}-m_i+\delta_{i,j})}}\sqrt{1-q^{2(m_{j+1}-m_j+1)}}
\ket{\underline{m}-\underline{\varepsilon}_i^j} \\
& ~ && \hskip-1cm\mathrm{if} \;0\leq i\leq j<k \;,\\
\pi^{(n)}_k(p_{ik})\ket{\underline{m}} &=q^{m_i+m_k}\sqrt{1-q^{2(m_{i+1}-m_i)}}
\ket{\underline{m}-\underline{\varepsilon}_i^k}
&& \hskip-1cm\mathrm{if} \;0\leq i<k \;,\\
\pi^{(n)}_k(p_{kk})\ket{\underline{m}} &=q^{2m_k}\ket{\underline{m}}
\;,\\
\pi^{(n)}_k(p_{ij})\ket{\underline{m}} &=0
&&\hskip-1cm\mathrm{if} \;j>k\;, 
\end{align*}
and $\pi^{(n)}_k(p_{ij})\ket{m}=0$ if $\underline{m}$ does not satisfy (\ref{eq:mconstr}).

Each representaton $\pi^{(n)}_k$ is an irreducible $*$-representation of both
$\A(S^{2n+1}_q)$ and $\A(\CP^n_q)$ when restricted to $\mathcal{V}^n_k$, and is identically zero outside $\mathcal{V}^n_k$.

\begin{lemma}\label{le:pl}
The spaces $\mathcal{V}^n_k$ enjoy the properties: 
$\mathcal{V}^n_j\perp\mathcal{V}^n_k$ if $|j-k|>1$, while
$\mathcal{V}^n_{k-1}\cap\mathcal{V}^n_k$ is the span of vectors $\ket{\underline{m}}$
with
\begin{equation}\label{eq:capconstr}
0\leq m_1\leq m_2\leq\ldots\leq m_k\;,\qquad m_k>m_{k+1}>\ldots>m_n\geq 0 \;,
\end{equation}
for all $1\leq k\leq n$.
As a consequence, if $|j-k|>1$, for all $a,b\in\A(S^{2n+1}_q)$
$$\pi^{(n)}_j(a) \,\pi^{(n)}_k(b)=0.$$ 
\end{lemma}

\begin{proof}
Let $0\leq j\leq k-2\leq n-2$. Due to (\ref{eq:mconstr}),
if $m_{k-1}\leq m_k$ the vector $\ket{\underline{m}}$ is not in $\mathcal{V}^n_j$,
and if $m_{k-1}>m_k$ the vector $\ket{\underline{m}}$ is not in $\mathcal{V}^n_k$.
This proves that $\mathcal{V}^n_j$ and $\mathcal{V}^n_k$ are orthogonal subspaces
of $\HH_n$. The remaining claims are straightforward.
\end{proof}

As a corollary of previous Lemma, the maps $\pi_\pm^{(n)}:\A(S^{2n+1}_q)\to\B(\HH_n)$, defined by
\begin{equation}\label{irreps+-}
\pi_+^{(n)}(a):=\sum_{\substack{0\leq k\leq n \\ k\;\mathrm{even}}}\pi^{(n)}_k(a) \;,\qquad
\pi_-^{(n)}(a):=\sum_{\substack{0\leq k\leq n \\ k\;\mathrm{odd}}}\pi^{(n)}_k(a) \;,
\end{equation}
are representations of the algebra $\A(S^{2n+1}_q)$.

\begin{prop}\label{prop}
The difference $\pi_+^{(n)}(a)-\pi_-^{(n)}(a)$ is of trace class on $\HH_n$ for all $a\in\A(\CP^n_q)$; furthermore, the trace is given by a series which -- as a function of $q$ -- is absolutely convergent on the open interval $0<q<1$.
\end{prop}

\begin{proof}
It is enough to prove the claim for $a=p_{ij}$, with $0\leq i\leq j\leq n$.
The space $\HH_n$ is the orthogonal direct sum of $\mathcal{V}^n_{k-1}\cap\mathcal{V}^n_k$, for all $1\leq k\leq n$, plus the joint kernel of all the representations involved.
By Lemma~\ref{le:pl}, on $\mathcal{V}^n_{k-1}\cap\mathcal{V}^n_k$ only the representations $\pi^{(n)}_{k-1}$ and $\pi^{(n)}_k$ are different from zero (and one contributes to $\pi_+^{(n)}$, while the other to $\pi_-^{(n)}$, according to the parity of $k$). We need to prove that $\pi^{(n)}_{k-1}(p_{ij})-\pi^{(n)}_k(p_{ij})$ is of trace class, and that the trace is absolutely convergent for any $0<q<1$.

{}From the explicit expressions given above, we see that $\pi^{(n)}_{k-1}(p_{ij})$ and $\pi^{(n)}_k(p_{ij})$ are both zero, if $j>k$. For $j=k$, $\pi^{(n)}_{k-1}(p_{ik})$ is zero and $\pi^{(n)}_k(p_{ik})$ has matrix coefficients
bounded by $q^{m_k}$. For $j=k-1$, one uses the inequality $|1-\sqrt{1-x^2}|\leq x$ --- valid for $0\leq x\leq 1$ --- to prove that $\pi^{(n)}_{k-1}(p_{ij})-\pi^{(n)}_k(p_{ij})$ still has matrix coefficients bounded by $q^{m_k}$. For $0\leq i\leq j\leq k-2$, the operators $\pi^{(n)}_{k-1}(p_{ij})$ and $\pi^{(n)}_k(p_{ij})$ coincide on $\mathcal{V}^n_{k-1}\cap\mathcal{V}^n_k$. The observation that the series
$$
\sum_{\underline{m}\;\textrm{satisfying (\ref{eq:capconstr})}}q^{m_k}
=\sum_{m_k=n-k}^\infty\binom{m_k+k-1}{k-1}\binom{m_k}{n-k}q^{m_k}
$$
is absolutely convergent for $0<q<1$ concludes the proof.
\end{proof}

Assembling things together, a $1$-summable even Fredholm module for $\A(\CP^n_q)$ is obtained
by using the representation $\pi^{(n)}:=\pi_+^{(n)}\oplus\pi_-^{(n)}$ on $\HH_{(n)} := \HH_n\oplus\HH_n$, the obvious grading operator $\gamma_{(n)}$
and, as usual,  
\begin{equation}\label{fmn}  
F_{(n)}:=\ma{14pt}{0 & 1 \\ 1 & 0} \;.
\end{equation}
All this is valid for any $n\geq 1$. Additional 
$n-1$ Fredholm modules are obtained by the same construction for $\A(\CP^j_q)$, $1\leq j\leq n-1$, and by pulling them back
to $\A(\CP^n_q)$. The last Fredholm module is the pullback of the canonical non-trivial Fredholm
module on $\C$, given on $\C\oplus\C$ by the representation $c\mapsto c\oplus 0$ and by
usual grading and operator $F$.

For $0\leq k\leq n$ and $N\in\N$, we denote $[\mu_k]$ the class of the Fredholm module 
$$
(\A(\CP^n_q), \, \HH_{(k)}, \, \pi^{(k)}, \, \gamma_{(k)}, \, F_{(k)}).
$$

% ----------------------------------------------------------------------
\section{K-theory}
We need some notation. The $q$-analogue of an integer number $n$ is  given by
$$
[n]:=\frac{q^n-q^{-n}}{q-q^{-1}} \;.
$$
This is defined for $q\not= 1$ and equals $n$ in the limit $q\to 1$. 
For any $n\geq 1$, define 
$$
[n]!:=[n][n-1]\ldots [1] \;,
$$
to be the $q$-factorial, with $[0]!:=1$, and the $q$-multinomial coefficients as
$$
[j_0,\ldots,j_n]!:=\frac{[j_0+\ldots+j_n]!}{[j_0]!\ldots[j_n]!} \;.
$$
With $N\geq 0$, let $\Psi_N=(\psi^N_{j_0,\ldots,j_n})$ be the vector-valued function on $S^{2n+1}_q$ with components
$$
\psi^N_{j_0,\ldots,j_n}:=[j_0,\ldots,j_n]!^{\frac{1}{2}}q^{-\frac{1}{2}\sum_{r<s}j_rj_s}
(z_n^{j_n}\ldots z_0^{j_0})^* \;,\qquad
\forall\;j_0+\ldots+j_n=N \;;
$$
there are $\binom{N+n}{n}$ of them. Then $\Psi_N^\dag\Psi_N=1$ and $P_N:=\Psi_N\Psi_N^\dag$ is a projection; the proof is in \cite{DD09}, and is a generalization of the case $n=2$ in \cite{DL09}. In particular $P_1=P$ is the `defining' projection in Sec.~\ref{se:intro} of the algebra $\A(\CP^n_q)$. 

Projections $P_{-N}$ are then obtained from $P_N$ with simple substitutions. 
We notice that the elements $z_i'=q^iz_i^*$ satisfy the defining relations of $\A(S^{2n+1}_{q^{-1}})$; 
the non-trivial relation to check is \eqref{eq:3rd}. Multiplying \eqref{eq:3rd} by $q^{2i}$ and summing for $i\geq k$ it becomes
$$
\sum\nolimits_{i\geq k}q^{2i}z_i^*z_i=q^{2k}\sum\nolimits_{i\geq k}z_iz_i^* \;,
$$
now subtracting $q^{-2}$ times the same equation evaluated for $k+1$, we get
$$
q^{2k}z_k^*z_k+(1-q^{-2})\sum\nolimits_{i>k}q^{2i}z_i^*z_i=q^{2k}z_kz_k^* \;,
$$
that means
$$
[z'^*_i,z'_i]=(1-q^{-2})\sum\nolimits_{j=i+1}^nz'_jz'^*_j \;.
$$

Vectors $\Psi_{-N}=(\psi^{-N}_{j_0,\ldots,j_n})$  are obtained out of 
$\Psi_N$ by first replacing $q$ with $q^{-1}$, and then replacing $z_i$ with $z'_i=q^iz_i^*$. Since the $q$-analogue is invariant under $q\to q^{-1}$, we obtain:
\begin{align*}
\psi^{-N}_{j_0,\ldots,j_n} &:=[j_0,\ldots,j_n]!^{\frac{1}{2}}q^{\frac{1}{2}\sum_{r<s}j_rj_s}
(z'^{j_0}_0)^*\ldots (z'^{j_n}_n)^* \\
&\:=[j_0,\ldots,j_n]!^{\frac{1}{2}}q^{\frac{1}{2}\sum_{r<s}j_rj_s+\sum_{r=0}^nrj_r}
z_0^{j_0}\ldots z_n^{j_n} \;,\qquad\forall\;j_0+\ldots+j_n=N \;.
\end{align*}
Then $\Psi_{-N}^\dag\Psi_{-N}=1$ and $P_{-N}:=\Psi_{-N}\Psi_{-N}^\dag$ is a projection.

For $N\in\N$ we denote by $[P_{-N}]$ be the class of the projection $P_{-N}$.
In the previous section we gave Fredholm module classes $[\mu_k]$ for each $0\leq k\leq n$.

\begin{prop}
For all $N\in\N$ and for all $0\leq k\leq n$ it holds that 
$$
\inner{[\mu_k], [P_{-N}]}:=
\tr_{\HH_k}(\pi_+^{(k)}-\pi_-^{(k)})(\tr\,P_{-N})
=\tbinom{N}{k} \;,
$$
with $\binom{N}{k}:=0$ when $k>N$.
\end{prop}
\begin{proof}
We have
$$
\inner{[\mu_k],[P_{-N}]}=\sum_{\underline{m}\in\N^k;m_k=0}\;\;\sum_{j=0}^k(-1)^j
\left<\underline{m}\right|\pi_j(\tr\,P_{-N})\ket{\underline{m}} \;.
$$
The integer $\inner{[\mu_0],[P_{-N}]}$ is the trace of the projection $P_{-N}$ evaluated at
the classical point (the unique $1$-dimensional representation):
one finds $\inner{[\mu_0],[P_{-N}]}=1$. The integer $\inner{[\mu_k],[P_{-N}]}$ being given,
for $1\leq k\leq n$, by a series which is absolutely convergent -- as a function of $q$ --
in the open interval $]0,1[$, it can be computed in the $q\to 0^+$ limit (cf. Sec.~5.3 in \cite{DD06} or \cite{DDLW07}). The case $N=0$ is trivial;  we focus on $N\geq 1$.
We notice that $q^{\frac{1}{2}n(n-1)}[n]!=1+O(q)$ for all $n\geq 0$. This implies
$$
q^{\sum_{r<s}j_rj_s}[j_0,\ldots,j_n]!=1+O(q)
$$
for all $j_0,\ldots,j_n\geq 0$. Thus,
$$
(\psi^{-N}_{j_0,\ldots,j_n})^*=\bigl\{1+O(q)\bigr\}\, q^{\sum_{r=0}^nrj_r}(z_n^{j_n})^*\ldots (z_0^{j_0})^* \;.
$$
Since matrix elements of $z_i, z_i^*$ (in the representations) are $O(1)$, we can conclude that
$(\psi^{-N}_{j_0,\ldots,j_n})^*=O(q)$ unless $\sum_{r=0}^nrj_r=0$,
which is equivalent to $j_0=N$. Hence, for the pairing with the $k$-th
Fredholm module, when $1\leq k\leq n$, we get, for any $0\leq j\leq k$, 
$$
\left<\underline{m}\right|\pi^{(k)}_j\bigl(\tr(P_{-N})\bigr)\ket{\underline{m}}
=\left<\underline{m}\right|\pi^{(k)}_j\bigl(z_0^N(z_0^N)^*\bigr)\bigr)\ket{\underline{m}}
+O(q) \;.
$$
In turn, 
\begin{multline*}
\inner{[\mu_k],[P_{-N}]}= \\
= \lim_{q\to 0^+}\sum_{j=1}^k(-1)^{j-1}\tr_{\mathcal{V}^k_{j-1}\cap\mathcal{V}^k_j}
\left\{
\left<\underline{m}\right|\pi^{(k)}_{j-1}\bigl(z_0^N(z_0^N)^*\bigr)\ket{\underline{m}}
-\left<\underline{m}\right|\pi^{(k)}_j\bigl(z_0^N(z_0^N)^*\bigr)\ket{\underline{m}}\right\} \;.
\end{multline*}
But $\pi^{(k)}_j\bigl(z_0^N(z_0^N)^*\bigr)$ is the identity operator on $\mathcal{V}^k_0$,
while for $j\geq 1$ we have
$$
\pi^{(k)}_j\bigl(z_0^N(z_0^N)^*\bigr)\ket{\underline{m}}=\begin{cases}
(1-q^{2m_1})(1-q^{2(m_1-1)})\ldots (1-q^{2(m_1-N+1)})\ket{\underline{m}}
& \mathrm{if}\;m_1\geq N \;, \\
0 & \mathrm{otherwise} \;.
\end{cases}
$$
Then, in the difference $\pi^{(k)}_{j-1}\bigl(z_0^N(z_0^N)^*\bigr)-\pi^{(k)}_j\bigl(z_0^N(z_0^N)^*\bigr)$,
the unique non-zero contribution comes from the term with $j=1$. Since
$$
\left<\underline{m}\right|\pi^{(k)}_1\bigl(z_0^N(z_0^N)^*\bigr)\ket{\underline{m}}=
\begin{cases}
1+O(q) & \mathrm{if}\;m_1\geq N \;, \\
0 & \mathrm{otherwise} \;,
\end{cases}
$$
we get
\begin{align*}
\inner{[\mu_k],[P_{-N}]} &=\lim_{q\to 0^+}\tr_{\mathcal{V}^k_0\cap\mathcal{V}^k_1}
\left\{1-\left<\underline{m}\right|\pi^{(k)}_1\bigl(z_0^N(z_0^N)^*\bigr)\ket{\underline{m}}\right\} \\
&=\sum_{N>m_1>m_2>\ldots>m_k\geq 0}1
=\binom{N}{k} \;,
\end{align*}
with the notation $\binom{N}{k}=0$ if $k>N$.
This concludes the proof.
\end{proof}

\begin{prop}
The elements $[\mu_0],\ldots,[\mu_n]$  are generators of $K^0(\A(\CP^n_q))$, and the elements
$[P_0],\ldots,[P_{-n}]$ are generators of $K_0(\A(\CP^n_q))$.
\end{prop}

\begin{proof}
Let $M\in \mathrm{Mat}_{n+1}(\Z)$ be the matrix having entries $M_{ij}:=\inner{[\mu_i],[P_{-j}]}$, for $i,j=0,1,\ldots,n$. Since
$$
\sum\nolimits_{k=j}^i(-1)^{j+k}\tbinom{i}{k}\tbinom{k}{j}=\begin{cases}
1 &\mathrm{if}\;i=j\;,\\
\sum_{r=i-k=0}^{i-j}(-1)^r\tbinom{i-j}{r}\tbinom{i}{j}=(1-1)^{i-j}\tbinom{i}{j}=0
&\mathrm{if}\;i>j\;,
\end{cases}
$$
the matrix $M$ has inverse $M^{-1}\in GL(n+1,\Z)$ with matrix entries $(M^{-1})_{ij}=(-1)^{i+j}\binom{j}{i}$.
This proves that the above-mentioned elements are a basis of $\Z^{n+1}$ as a $\Z$-module, which is equivalent to say that they generate $\Z^{n+1}$ as abelian group.
\end{proof}

% ----------------------------------------------------------------------

\section{On Dirac operators}

For any $a\in\A(\CP^n_q)$, the matrix coefficients of $\pi_+^{(n)}(a)-\pi_-^{(n)}(a)$, 
where $\pi_\pm$ are the representations in \eqref{irreps+-}, go to
zero exponentially. This fact suggests using these representations to
construct even spectral triples for $\A(\CP^n_q)$
(as was done for the standard Podle\'s sphere in Prop.~5.1 in \cite{DDLW07})
with Dirac operator having sign given by \eqref{fmn}.
One obtains spectral triples of any summability $d\in\R^+$. 
For example, an $n$-summable spectral triple is obtained by defining the operator 
$D=|D|F_{(n)}$ on $\HH_{(n)} = \HH_n\oplus\HH_n$ as
$$
|D|\ket{\underline{m}}:=(m_1+\ldots+m_n)\ket{\underline{m}} \;.
$$
The multiplicity of the eigenvalue $\pm\lambda$, $\lambda\in\N$, is $\binom{\lambda+n}{n-1}$. This is a polynomial in $\lambda$ of order $n-1$, and so the metric dimension is $n$, as claimed.
The bounded commutators condition follows from the identity
$
[D,a]=|D|[F_{(n)},a]+[|D|,a]F_{(n)}
$
and the observation that:  
\begin{itemize}
\item[i)] the first term is bounded, since $|D|$ has a spectrum
that is polynomially divergent and $[F_{(n)},a]$ has matrix coefficients that
go to zero exponentially,  
\item[ii)] the second term is bounded as well, since one sees that 
the explicit expression of $\pi_\pm(p_{ij})$ are sums of bounded shift operators, which are eigenvectors of $[|D|,\,\,\cdot\,\,]$.
\end{itemize}

With some extra work, by setting 
$|D|\ket{\underline{m}}:=(m_1+\ldots+m_n)^{n/d}\ket{\underline{m}}$ for any $d\in\R^+$, 
one gets a $d$-summable spectral triple. We expect all these $d$-summable triples to be regular. 

On the other hand, $0^+$-summable (non-regular) spectral triples
for quantum projective spaces have been recently constructed in \cite{DD09}. 

% ----------------------------------------------------------------------

\bigskip\medskip
\begin{center}
\textsc{Acknowledgments}.
\end{center}
GL was partially supported by the `Italian project Cofin06 - Noncommutative geometry, quantum groups and applications'.

\bigskip
\providecommand{\bysame}{\leavevmode\hbox to3em{\hrulefill}\thinspace}

\end{document}